\documentclass[11pt,oneside]{amsart}

\usepackage{amssymb,amsthm,amsmath}
\usepackage{pstricks}
\usepackage{pst-node}
\usepackage{pst-tree}
\usepackage{graphicx}
\usepackage{pstcol}
\usepackage{pst-plot}
\usepackage{subfigure}
\usepackage{color}
\definecolor{darkblue}{rgb}{0, 0, .4}
\usepackage[bookmarks,breaklinks]{hyperref}
\hypersetup{
	colorlinks=true,
	linkcolor=darkblue,
	anchorcolor=darkblue,
	citecolor=darkblue,
	pagecolor=darkblue,
	urlcolor=darkblue,
	pdfpagemode=UseThumbs,
	pdftitle={Finitely labeled generating trees and restricted permutations},
	pdfsubject={Enumerative Combinatorics},
	pdfauthor={Vincent Vatter},
	pdfkeywords={restricted permutation, forbidden subsequence, generating tree}
}
	
\newtheorem{theorem}{Theorem}
\newtheorem{proposition}[theorem]{Proposition}
\newtheorem{lemma}[theorem]{Lemma}

\newtheorem{remark}[theorem]{Remark}

\newcounter{todocounter}

\newcommand{\st}{\operatorname{st}}
\newcommand{\Av}{\operatorname{Av}}
\newcommand{\R}{\mathcal{R}}
\newcommand{\Binf}{\|B\|_\infty}

\begin{document}

\title{Finitely labeled generating trees and restricted permutations}

\author[Vincent Vatter]{Vincent Vatter}
\address{School of Mathematics and Statistics, University of St.~Andrews, St.~Andrews, Fife, Scotland}
\email{\href{mailto:vince@mcs.st-and.ac.uk}{\texttt{vince@mcs.st-and.ac.uk}}}
\urladdr{\url{http://www-groups.mcs.st-and.ac.uk/~vince/}}
\thanks{Partially supported by an award from DIMACS and an NSF VIGRE grant to the Rutgers University Department of Mathematics.}

\date{\today}

\begin{abstract}
Generating trees are a useful technique in the enumeration of various combinatorial objects, particularly restricted permutations.  Quite often the generating tree for the set of permutations avoiding a set of patterns requires infinitely many labels.  Sometimes, however, this generating tree needs only finitely many labels.  We characterize the finite sets of patterns for which this phenomenon occurs.  We also present an algorithm --- in fact, a special case of an algorithm of Zeilberger --- that is guaranteed to find such a generating tree if it exists.
\end{abstract}

\maketitle

\section{Introduction}\label{finlabel-intro}

The {\it standardization\/} of a word $w$ consisting of $k$ distinct integers is the permutation $\st(w)$ of length $k$ obtained by replacing the smallest entry of $w$ by $1$, the second smallest entry by $2$, and so on.  If $\beta\in S_k$, we say that the length $n$ word of distinct integers $w$ {\it contains a $\beta$-pattern\/} if and only if it contains a (not necessarily contiguous) subword whose standardization is $\beta$.  Otherwise we say that $w$ is $\beta$-avoiding.  For example, the permutation $239145678$ contains a $312$-pattern (for example, the subword $916$ standardizes to $312$) but avoids $321$.  If $B$ is a set of permutations, we say that the word $w$ avoids $B$ if it avoids every member of $B$.  Such a set $B$ is often referred to as a set of {\it forbidden patterns\/}.  Although we have allowed $w$ to be a word in our definition, we are most interested in the case where $w$ is a permutation.  Let $\Av_n(B)$ denote the set of $B$-avoiding permutations of length $n$ and let $\Av(B)$ denote the set of all finite $B$-avoiding permutations.  We refer to $\sum_n |\Av_n(B)|x^n$ as the generating function for $\Av(B)$.

The problem of enumerating restricted permutations has received considerable attention over the last two decades.  B\'ona's book~\cite{bona:book} provides an overview of these efforts.  A common technique is that of generating trees, introduced by Chung, Graham, Hoggatt, and Kleiman~\cite{cghk:baxter} and used extensively by many others\footnote{To name a few: Barcucci, Del Lungo, Pergola, and Pinzani~\cite{bdpp:motz, bdpp:incnum}; Chen, Mansour, and Yan~\cite{cmy:map}; Chow and West~\cite{cw:cheby}; Dulucq, Gire, and Guibert~\cite{dgg:west}; Dulucq, Gire, and West~\cite{dgw:maps}; Guibert and Pergola~\cite{gp:vex}; Kremer~\cite{k:sn}; Kremer and Shiu~\cite{ks:len4}; Marinov and Radoi{\v{c}}i{\'c}~\cite{mr:1324}; Merlini, Sprugnoli, and Verri~\cite{msv:tennis}; Pergola and Sulanke~\cite{ps:schroder}; Stankova~\cite{s:len4,stankova:fs}; Stankova and West~\cite{sw:321hex}; and West~\cite{west:cat,west:trees}.} since.  The recently introduced ECO (enumerating combinatorial objects) method\footnote{See Barcucci, Del Lungo, Pergola, and Pinzani~\cite{eco:survey} for a survey, and Barcucci, Pergola, Pinzani, and Rinaldi~\cite{eco:hills} or Ferrari, Pinzani, and Rinaldi~\cite{eco:parts} for other applications.} extends the notion of generating trees to other combinatorial contexts.  There has also been some interest in the algebraic properties of generating trees and ECO systems\footnote{For this the reader is referred to Duchi, Fedou, and Rinaldi~\cite{eco:grammars}; Ferrari, Pergola, Pinzani, and Rinaldi~\cite{eco:alg}; and Merlini, Sprugnoli, and Verri~\cite{msv:alg-gt}.}.

Precisely, a generating tree is a rooted, labeled, and typically infinite tree such that the label of a node determines the labels of its children.  Sometimes the labels of the tree are taken to be natural numbers, but this is not necessary, and the algorithm we will describe labels nodes by permutations.  Therefore we specify a generating tree by supplying the label of the root and a set of {\it succession rules\/}.  For example, the complete binary tree may be given by
$$
\begin{array}{llcl}
\mbox{Root:}	&	(SS)&&\\
\mbox{Rule:}	&	(SS)&\leadsto &(SS)(SS).
\end{array}
$$

\begin{figure}[t]
\begin{footnotesize}
\begin{center}
\psset{xunit=0.028in, yunit=0.02in}
\psset{linewidth=0.25\psxunit}
\begin{pspicture}(-5,5)(208,78)
\pscircle*(10,10){1\psxunit}
\pscircle*(25,10){1\psxunit}
\pscircle*(40,10){1\psxunit}
\pscircle*(55,10){1\psxunit}
\pscircle*(70,10){1\psxunit}
\pscircle*(85,10){1\psxunit}
\pscircle*(100,10){1\psxunit}
\pscircle*(115,10){1\psxunit}
\pscircle*(130,10){1\psxunit}
\pscircle*(145,10){1\psxunit}
\pscircle*(160,10){1\psxunit}
\pscircle*(175,10){1\psxunit}
\pscircle*(190,10){1\psxunit}
\rput[c](10,5){$1234$}
\rput[c](25,5){$4123$}
\rput[c](40,5){$3124$}
\rput[c](55,5){$3412$}
\rput[c](70,5){$4312$}
\rput[c](85,5){$2134$}
\rput[c](100,5){$4213$}
\rput[c](115,5){$2314$}
\rput[c](130,5){$2341$}
\rput[c](145,5){$4231$}
\rput[c](160,5){$3214$}
\rput[c](175,5){$3421$}
\rput[c](190,5){$4321$}
\pscircle*(17.5,30){1\psxunit}
\pscircle*(55,30){1\psxunit}
\pscircle*(92.5,30){1\psxunit}
\pscircle*(130,30){1\psxunit}
\pscircle*(175,30){1\psxunit}
\rput[r](14.5,30){$123$}
\rput[r](52,30){$312$}
\rput[r](89.5,30){$213$}
\rput[l](133,30){$231$}
\rput[l](178,30){$321$}
\psline(10,10)(17.5,30)
\psline(25,10)(17.5,30)
\psline(40,10)(55,30)
\psline(55,10)(55,30)
\psline(70,10)(55,30)
\psline(85,10)(92.5,30)
\psline(100,10)(92.5,30)
\psline(115,10)(130,30)
\psline(130,10)(130,30)
\psline(145,10)(130,30)
\psline(160,10)(175,30)
\psline(175,10)(175,30)
\psline(190,10)(175,30)
\pscircle*(36.25,50){1\psxunit}
\pscircle*(130,50){1\psxunit}
\rput[c](30.25,54){$12$}
\rput[c](136,54){$21$}
\psline(17.5,30)(36.25,50)
\psline(55,30)(36.25,50)
\psline(92.5,30)(130,50)
\psline(130,30)(130,50)
\psline(175,30)(130,50)
\pscircle*(83.125,70){1\psxunit}
\rput[c](83.125,76){$1$}
\psline(36.25,50)(83.125,70)
\psline(130,50)(83.125,70)
\end{pspicture}
\end{center}
\end{footnotesize}
\caption{The first four levels of the pattern-avoidance tree $T(132,3241)$}\label{F-132-3241}
\end{figure}
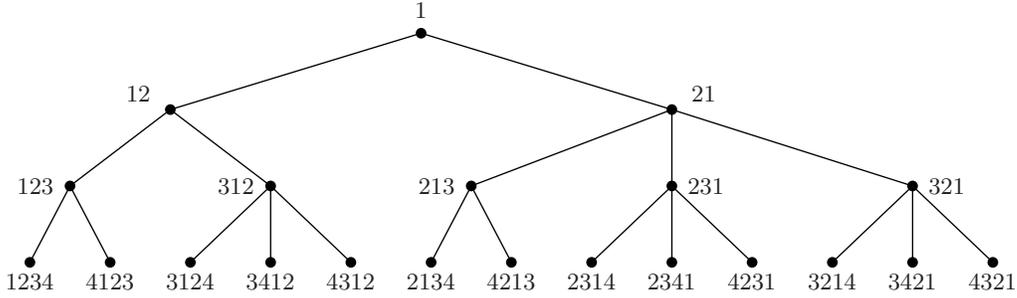

The connection to restricted permutations comes through {\it pattern-avoidance trees\/}.  We say that the permutation $\sigma$ of length $n$ is a {\it child\/} of $\pi\in S_{n-1}$ if $\sigma$ can be obtained by inserting $n$ into $\pi$.  This defines a rooted tree $T$ on the set of all permutations.  For a set of patterns $B$, we define the pattern-avoidance tree $T(B)$ to be the subtree of $T$ whose nodes are the $B$-avoiding permutations $\Av(B)$.  For example, the first four levels of $T(132,3241)$ are shown in Figure~\ref{F-132-3241}.  The {\it active sites\/} (relative to $B$) of $\pi\in \Av_{n-1}(B)$ are the positions $i$ for which inserting $n$ right before the $i$th entry of $\pi$ produces a $B$-avoiding permutation.  By convention, $n+1$ is an active site of $\pi$ if appending $n$ to the end of $\pi$ produces a $B$-avoiding permutation.  An inactive site is any site that is not active.  For example, the active sites of $213$ relative to $B=\{132,3241\}$ are $1$ and $4$, whereas the inactive sites are $2$ and $3$.

Given a pattern-avoidance tree $T(B)$, we would like to find an isomorphic (in the category of rooted trees) generating tree.  For example, take a permutation $\pi\in\Av_n(132,3241)$.  There are at most three sites in which we may insert $n+1$ to form a $\{132,3241\}$-avoiding child: the beginning, the end, and the site directly to the right of $n$.  Indeed, if we insert $n+1$ to the left of $n$ but not at the beginning we form a $132$-pattern, while if we insert $n+1$ further than one site to the right of $n$ but not at the end we get a subsequence $n, x, (n+1), y$ and either $x<y$, giving a $132$-pattern, or $x>y$, giving a $3241$-pattern.  We can therefore constuct an isomorphic generating tree with just two labels, label $(2)$ for the nodes where $n$ is the last entry and label $(3)$ for the nodes where $n$ is not the last entry.  The generating tree is then
$$
\begin{array}{llcl}
\mbox{Root:}	&	(2)&&\\
\mbox{Rules:}	&	(2)&\leadsto &(2)(3)\\
&(3)&\leadsto&(2)(3)(3).
\end{array}
$$

On the other hand, the pattern-avoidance tree $T(123)$ requires infinitely many labels; by considering the lexicographically first ascent in a $123$-avoiding permutation one can show that $T(123)$ is isomorphic to the generating tree given by
$$
\begin{array}{llcl}
\mbox{Root:}	&	(2)&&\\
\mbox{Rules:}	&	(j)&\leadsto &(2)(3)\cdots(j+1).
\end{array}
$$

Upon finding a generating tree isomorphic to $T(B)$, one often wishes to get the generating function for $\Av(B)$.  In general, as witnessed by Bousquet-M{\'e}lou~\cite{bm:four} and Banderier, Bousquet-M{\'e}lou, Denise, Flajolet, Gardy, and Gouyou-Beauchamps~\cite{gffgt}, this process can be quite intricate.  However, if the generating tree has only finitely many labels then the generating function for $\Av(B)$, which must be rational, can be computed using the transfer matrix method (see Stanley's text~\cite[Section 4.7]{stanley:ec1} for details).

Herein we characterize the finite sets $B$ for which $T(B)$ is isomorphic to a finitely labeled generating tree, answering a question raised earlier~\cite{2len3}.  One requirement for $T(B)$ to be isomorphic to a finitely labeled generating tree is that there must be a bound on the number of children a node of $T(B)$ may have.  For this to occur, $B$ must contain both a child of an increasing permutation (such as $132$, $4123$, or $12345$) and a child of a decreasing permutation (such as $231$, $3241$, or $54321$), because otherwise either $12\cdots n$ or $n\cdots 21$ will have $n+1$ children for all $n$.  In fact, Kremer and Shiu~\cite{ks:len4} showed that this is enough.  We include a short proof below.  First recall the following famous theorem of Erd\H os and Szekeres.

\begin{theorem}[Erd\H{o}s and Szekeres~\cite{es:acpig}]\label{es}
Every permutation of length at least $(k-1)(\ell-1)+1$ contains either an increasing sequence of length $k$ or a decreasing sequence of length $\ell$.
\end{theorem}

\begin{theorem}[Kremer and Shiu~\cite{ks:len4}]\label{bounded}
The pattern-avoidance tree $T(B)$ has bounded degrees if and only if $B$ contains both a child of an increasing permutation and a child of a decreasing permutation.
\end{theorem}
\begin{proof}
We have already noted that the condition on $B$ is necessary, so it suffices to show that it is sufficient.  Assume not, so that although $B$ contains a child of $12\cdots k$ and a child of $\ell\cdots 21$, $T(B)$ does not have bounded degrees.  Set $n=(k-1)(\ell-1)+1$.

We claim that there is a permutation $\pi\in T(B)$ of length $n$ with $n+1$ children (or, equivalently, $n+1$ active sites).  Since $T(B)$ does not have bounded degrees, we can find a permutation in $T(B)$ with at least $n+1$ active sites.  Suppose that $i_1<i_2<\cdots<i_{n+1}$ are active sites of $\pi$.  Now form the word $w=\pi({i_1})\pi({i_2})\cdots \pi({i_n})$ and set $\sigma=\st(w)$.  For example, suppose that $n=3$ and $\pi=461523$ with active sites $\{1,2,5,7\}$.  Then we get $w=\pi(1)\pi(2)\pi(5)=462$ and $\sigma=\st(462)=231$.

By construction, $\sigma$ is a permutation of length $n$ with $n+1$ children in $T(B)$, proving the claim.  However, by our choice of $n$, the Erd\H{o}s-Szekeres theorem shows that $\sigma$ contains either an increasing subsequence of length at least $k$ or a decreasing subsequence of length at least $\ell$.  Thus we have reached a contradiction because our assumptions on $B$ imply that at least one of the children of $\sigma$ most contain a pattern from $B$ and thus cannot be a node of $T(B)$.
\end{proof}

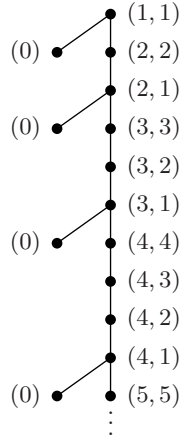
\begin{figure}[t]
\begin{footnotesize}
\begin{center}
\psset{xunit=0.028in, yunit=0.02in}
\psset{linewidth=0.25\psxunit}
\begin{pspicture}(-3,-7)(30,103)
\pscircle*(10,0){1\psxunit}
\pscircle*(10,10){1\psxunit}
\pscircle*(10,20){1\psxunit}
\pscircle*(10,30){1\psxunit}
\pscircle*(10,40){1\psxunit}
\pscircle*(10,50){1\psxunit}
\pscircle*(10,60){1\psxunit}
\pscircle*(10,70){1\psxunit}
\pscircle*(10,80){1\psxunit}
\pscircle*(10,90){1\psxunit}
\pscircle*(10,100){1\psxunit}
\rput[l](13,0){$(5,5)$}
\rput[l](13,10){$(4,1)$}
\rput[l](13,20){$(4,2)$}
\rput[l](13,30){$(4,3)$}
\rput[l](13,40){$(4,4)$}
\rput[l](13,50){$(3,1)$}
\rput[l](13,60){$(3,2)$}
\rput[l](13,70){$(3,3)$}
\rput[l](13,80){$(2,1)$}
\rput[l](13,90){$(2,2)$}
\rput[l](13,100){$(1,1)$}
\psline(10,0)(10,100)
\pscircle*(0,0){1\psxunit}
\pscircle*(0,40){1\psxunit}
\pscircle*(0,70){1\psxunit}
\pscircle*(0,90){1\psxunit}
\psline(0,0)(10,10)
\psline(0,40)(10,50)
\psline(0,70)(10,80)
\psline(0,90)(10,100)
\rput[r](-3,0){$(0)$}
\rput[r](-3,40){$(0)$}
\rput[r](-3,70){$(0)$}
\rput[r](-3,90){$(0)$}
\rput[c](10,-5){$\vdots$}
\end{pspicture}
\end{center}
\end{footnotesize}
\caption{A tree with bounded degrees that is not isomorphic to a finitely labeled generating tree}\label{F-bounded}
\end{figure}

In general, trees with bounded degrees need not be isomorphic to finitely labeled generating trees.  For example, consider the generating tree (pictured in Figure~\ref{F-bounded}) given by
$$
\begin{array}{llcll}
\mbox{Root:}	&	(1,1)&&\\
\mbox{Rules:}	&	(i,j)&\leadsto &(i,j-1)\mbox{ if $j\ge 2$,}\\
			&	(i,1)&\leadsto &(i+1,i+1)(0),
\end{array}
$$
where nodes labeled by $(0)$ do not produce children.  This tree is clearly not isomorphic to a finitely labeled generating tree since the distance between two nodes of degree 2 is unbounded, but each of its nodes has at most two children.  Our main result, Theorem~\ref{finlabel}, says that if $B$ is finite and $T(B)$ has bounded degrees then $T(B)$ is isomorphic to a finitely labeled generating tree, so this example shows that our proof will need to make use of the special properties of pattern-avoidance trees.

In the next section we introduce lemmas and notation for pattern-avoidance trees.  Section~\ref{finlabel-alg} describes our labeling algorithm, while the proof that this algorithm works is contained in Section~\ref{finlabel-proof}.

\section{Removable \& GT-reducible entries}\label{finlabel-lemmas}

In order to motivate our technique we begin by returning to the example of $T(132,3241)$.  Since our approach in the last section was rather ad hoc, we now attempt to analyze this tree (or rather, two of its nodes) in a more systematic manner.  Consider $213$.  In this permutation we can insert new maximal entries at the beginning or the end, but we will never be able to insert a new maximal entry between the $2$ and the $3$.  Now let $\pi$ denote a descendant of $213$.  From our previous comments, $213$ appears as a contiguous block in $\pi$.  Observe that for any possible $132$ or $3142$-pattern in $\pi$ involving the $1$ there is another pattern which uses the $2$ instead.  Therefore it does not matter whether or not the $1$ is stuck between the $2$ and the $3$, and we can assign the same label to $213$ as we assign to $12$.  In fact, a similar argument shows that $12$ can be labeled by the same label as $1$ receives.  It is notions like these that we aim to formalize in this section.

First, we say that an entry $x$ in the permutation $\pi$ is {\it removable\/} (relative to a set of patterns $B$) if it is adjacent to at most one active site.  For example, every entry of $213$ is removable when $B=\{132,3241\}$, although no entry of $21$ is removable.  When more detail is needed, we say that the entry $\pi(i)$ is left-removable if $i$ is an inactive site of $\pi$ and right-removable if $i+1$ is inactive, so every removable entry is either left-removable, right-removable, or both.

If $x$ is an entry of the word $w$ and $w$ contains distinct integers then we write $w-x$ to denote the word formed from $w$ by removing $x$.  For example, $461523-2=46153$.  If $X$ is a set of entries of $w$, we similarly write $w-X$ to denote the word formed by removing each entry in $X$.

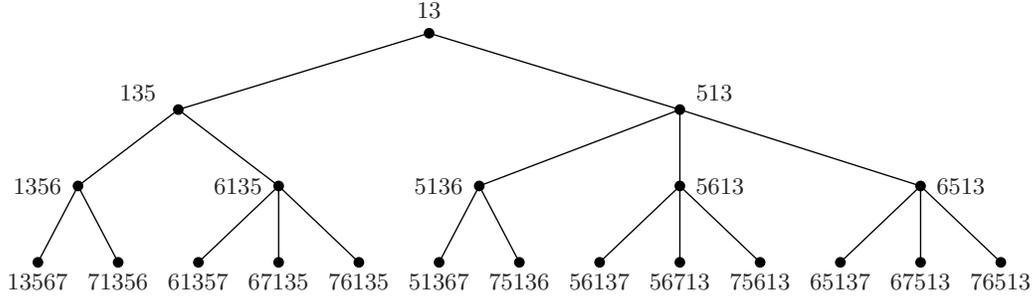
\begin{figure}[t]
\begin{footnotesize}
\begin{center}
\psset{xunit=0.028in, yunit=0.02in}
\psset{linewidth=0.25\psxunit}
\begin{pspicture}(-5,5)(208,78)
\pscircle*(10,10){1\psxunit}
\pscircle*(25,10){1\psxunit}
\pscircle*(40,10){1\psxunit}
\pscircle*(55,10){1\psxunit}
\pscircle*(70,10){1\psxunit}
\pscircle*(85,10){1\psxunit}
\pscircle*(100,10){1\psxunit}
\pscircle*(115,10){1\psxunit}
\pscircle*(130,10){1\psxunit}
\pscircle*(145,10){1\psxunit}
\pscircle*(160,10){1\psxunit}
\pscircle*(175,10){1\psxunit}
\pscircle*(190,10){1\psxunit}
\rput[c](10,5){$13567$}
\rput[c](25,5){$71356$}
\rput[c](40,5){$61357$}
\rput[c](55,5){$67135$}
\rput[c](70,5){$76135$}
\rput[c](85,5){$51367$}
\rput[c](100,5){$75136$}
\rput[c](115,5){$56137$}
\rput[c](130,5){$56713$}
\rput[c](145,5){$75613$}
\rput[c](160,5){$65137$}
\rput[c](175,5){$67513$}
\rput[c](190,5){$76513$}
\pscircle*(17.5,30){1\psxunit}
\pscircle*(55,30){1\psxunit}
\pscircle*(92.5,30){1\psxunit}
\pscircle*(130,30){1\psxunit}
\pscircle*(175,30){1\psxunit}
\rput[r](14.5,30){$1356$}
\rput[r](52,30){$6135$}
\rput[r](89.5,30){$5136$}
\rput[l](133,30){$5613$}
\rput[l](178,30){$6513$}
\psline(10,10)(17.5,30)
\psline(25,10)(17.5,30)
\psline(40,10)(55,30)
\psline(55,10)(55,30)
\psline(70,10)(55,30)
\psline(85,10)(92.5,30)
\psline(100,10)(92.5,30)
\psline(115,10)(130,30)
\psline(130,10)(130,30)
\psline(145,10)(130,30)
\psline(160,10)(175,30)
\psline(175,10)(175,30)
\psline(190,10)(175,30)
\pscircle*(36.25,50){1\psxunit}
\pscircle*(130,50){1\psxunit}
\rput[c](28.75,54.5){$135$}
\rput[c](136.5,54.5){$513$}
\psline(17.5,30)(36.25,50)
\psline(55,30)(36.25,50)
\psline(92.5,30)(130,50)
\psline(130,30)(130,50)
\psline(175,30)(130,50)
\pscircle*(83.125,70){1\psxunit}
\rput[c](83.125,76){$13$}
\psline(36.25,50)(83.125,70)
\psline(130,50)(83.125,70)
\end{pspicture}
\end{center}
\end{footnotesize}
\caption{The first four levels of $W(\{132,3241\};4,13)$}\label{F-132-3241-W}
\end{figure}

When $\pi$ is a node of $T(B)$ we let $T(B;\pi)$ denote the subtree consisting of $\pi$ and its descendants.  In order to avoid the shifting of indices and values caused by standardization we will also make use of the similar but less natural tree $W(B;n,u)$, which we define whenever $u$ is a $B$-avoiding word containing distinct integers all at most $n$.  The root of $W(B;n,u)$ is $u$, and it contains all $B$-avoiding words $w$ that can be formed by shuffling $u$ with a permutation on $[n+1,n+|w|-|u|]$.  If $v$ and $w$ are nodes of $W(B;n,u)$ and $w$ has greatest entry $m$, then $w$ is a child of $v$ if $w$ can be obtained by inserting $m$ into $v$, that is, if $v=w-m$.  An example is shown in Figure~\ref{F-132-3241-W}.  For a word $w\in W(B;n,u)$ we define active sites, inactive sites, removability, left-removability, and right-removability as we did for the permutations of $T(B)$.

Once these definitions are unraveled, it is evident that $W(B;n,u)\cong T(B;\st(u))$ for all allowed values of $n$, and if $\pi$ is a permutation then $W(B;|\pi|,\pi)$ and $T(B;\pi)$ are not merely isomorphic, but they are the same tree.

If $u$ is a word of length at least two containing the entry $x$ we define the map
$$
\partial_{x}:W(B;n,u)\rightarrow W(B;n,u-x)
$$
by $\partial_{x}(w)=w-x$.  Our upcoming Proposition~\ref{embed} shows that $\partial_x$ is one-to-one if $x$ is a removable entry in $u$.  There are several different cases in the definition of the inverse map.

First suppose that $x=u(i)$ for some $i>2$.  In this case we define the map
$$
\iota^-_{u,x}:W(\emptyset;n,u-x)\rightarrow W(\emptyset;n,u)
$$
by letting $\iota^-_{u,x}(w)$ denote the word obtained from $w$ by inserting $x$ immediately to the right of the entry $u(i-1)$.  If $x$ is left-removable in $u$ relative to $B$ then it is easy to see that $\iota^-_{u,x}$ maps words in $\partial_x(W(B;n,u))$ to words in $W(B;n,u)$.  Furthermore, if $w\in W(B;n,u-x)\setminus \partial_x(W(B;n,u))$ then $\iota^-_{u,x}(w)$ will contain at least one pattern from $B$.

There are three more cases for us to define this operation.  If $x=u(1)$ then we simply let $\iota^-_{u,x}(w)=xw$.  Similarly, if $x$ is the last entry of $u$ we let $\iota^+_{u,x}(w)=wx$.  Otherwise we let $\iota^+_{u,x}(w)$ denote the word formed from $w$ by inserting $x$ to the immediate left of $u(i+1)$.
\begin{proposition}\label{embed}
Let $u$ be a $B$-avoiding word of length at least two containing distinct integers all at most $n$ and suppose that $x\in u$ is removable.  Then  $\partial_{x}:W(B;n,u)\rightarrow W(B;n,u-x)$ is one-to-one.  More specifically, if $x$ is left-removable then
$$
\iota^-_{u,x}\circ\partial_{x}:W(B;n,u)\rightarrow W(B;n,u)
$$
is the identity, and if $x$ is right-removable then
$$
\iota^+_{u,x}\circ\partial_{x}:W(B;n,u)\rightarrow W(B;n,u)
$$
is the identity.
\end{proposition}
\begin{proof}
As the various cases are quite similar, let us assume that $x=u(i)$ for some $i>2$ and that $x$ is left-removable.  Take $w\in W(B;n,u)$.  Since $x$ is left-removable, $w$ cannot contain an entry between $u(i-1)$ and $x$.  Thus applying $\partial_x$ removes $x$, but then $\iota^-_{u,x}$ inserts $x$ immediately to the right of $u(i-1)$, restoring $w$.
\end{proof}

One of the implications of Proposition~\ref{embed} is that $W(B;n,u)$ embeds into $W(B;n,u-x)$ by the map $\partial_x$ whenever $x$ is a removable entry of $u$.  It sometimes happens that these two trees are isomorphic.  If $W(B;n,u)\cong W(B;n,u-x)$, and $x$ is removable, then we say that $x$ is {\it generating-tree-reducible relative to $B$\/}, or for short, {\it GT-reducible\/}.  This is a strengthening of Zeilberger's definition of reversely deleteable entries from \cite{z:wilf}.  In this language, our observation at the beginning of the section was that the entry $1$ in $213$ is GT-reducible for $B=\{132,3241\}$.  If the permutation $\pi$ contains a GT-reducible entry, we also refer to $\pi$ as being GT-reducible.

Before ending our discussion of $\iota$, let us note that in many cases these maps commute:
\begin{remark}\label{commute}
Let $u$ be a word containing distinct integers all at most $n$, let $x$ and $y$ be nonadjacent entries of $u$, and let $\delta,\epsilon\in\{+,-\}$.  Then
$$
\iota^\delta_{u,x}\circ\iota^\epsilon_{u-x,y}=\iota^\epsilon_{u,y}\circ\iota^\delta_{u-y,x}
$$
as maps from $W(\emptyset;n,u-x-y)$ to $W(\emptyset;n,u)$.
\end{remark}

We would now like to show that it is possible to decide whether a removable entry is GT-reducible.  For this we need two more definitions.  Given a tree $T$, let $\ell_i(T)$ denote the number of nodes of $T$ of height $i$, so $\ell_0(T)=1$ unless $T$ is the empty tree, $\ell_1(T)$ is the number of children of the root node, and $\ell_2(T)$ is the number of grandchildren of the root node.  Also, if $B$ is a finite set of patterns, let $\Binf$ denote the length of the longest pattern in $B$.

\begin{proposition}\label{reinsert}
Let $u$ be a $B$-avoiding word of length at least two containing distinct integers all at most $n$.  The removable entry $x\in u$ is GT-reducible if and only if
$$
\ell_r(W(B;n,u))=\ell_r(W(B;n,u-x))
$$
for all $1\le r\le \Binf-1$.
\end{proposition}
\begin{proof}
If $x\in u$ is GT-reducible then $W(B;n,u)\cong W(B;n,u-x)$ by definition, so $\ell_r(W(B;n,u)) = \ell_r(W(B;n,u-x))$ for all $r\in\mathbb{N}$.  Suppose now that $\ell_r(W(B;n,u)) = \ell_r(W(B;n,u-x))$ for all $1\le r\le \|B\|_{\infty}-1$.  Since $x$ is removable, Proposition~\ref{embed} shows that $W(B;n,u)$ embeds into $W(B;n,u-x)$ by the map $\partial_{x}$.  Thus we would like to show that this map is onto.

Suppose not and choose $w\in W(B;n,u-x)$ that is not in the image of $\partial_{x}$.  Let $\epsilon=-$ if $x$ is left-removable in $u$.  Otherwise $x$ must be right-removable in $u$, and here we let $\epsilon=+$.  Since $w\notin \partial_x(W(B;n,u))$, $\iota^\epsilon_{u,x}(w)$ contains a permutation from $B$.  Choose a subword of $\iota^\epsilon_{u,x}(w)$ that standardizes to a member of $B$ and label it $v$.  Because $w$ avoids $B$, $x$ must be an entry of $v$.

Now consider the subword $w'$ of $\iota^\epsilon_{u,x}(w)$ containing all entries that are either in $u$ or in $v$.  Because $w'$ contains $v$, it contains a pattern from $B$.  However, $w'-x$ is a subword of $w$, so it avoids $B$.  We would like to find a word in $\iota^\epsilon_{u,x}(W(B;n,u-x))$ with these properties.  We do this by ``partially standardizing'' $w'$: replace the smallest entry of $w'$ that is not in $u$ by $n+1$, the next smallest entry of $w'$ that is not in $u$ by $n+2$, and so on.  Label the resulting word $w''$.  Notice that $w''\in \iota^\epsilon_{u,x}(W(B;n,u-x))$, $w''$ contains a pattern from $B$, and $w''-x$ is $B$-avoiding.  These observations and Proposition~\ref{embed} show that
$$
\ell_{|w''|-|u|}(W(B;n,u))<\ell_{|w''|-|u|}(W(B;n,u-x)).
$$
Furthermore, since both $w''$ and $u$ contain $x$, $|w''|\le |u|+\Binf-1$, contradicing our hypotheses.
\end{proof}

\section{The algorithm}\label{finlabel-alg}

Our work in the previous section suggests the following approach for finding a generating tree isomorphic to $T(B)$.  If $1\notin \Av(B)$ then our task is quite easily accomplished, so let us assume that $1$ avoids $B$.  We start with a root node $(1)$, a set $P=\{1\}$ of permutations that we have not checked for GT-reducible entries, and a set $\R=\emptyset$ of succession rules.  Now we pick a permutation $\pi\in P$ of minimum length and check it for GT-reducible entries (we make the convention that the permutation $1$ never has a GT-reducible entry).

First suppose that $\pi$ is not GT-reducible, and say that its $B$-avoiding children are $\sigma_1,\sigma_2,\dots,\sigma_t$.  In this case we remove $\pi$ from $P$, add its $B$-avoiding children to $P$, and add the succession rule
$$
(\pi)\leadsto(\sigma_1)(\sigma_2)\cdots(\sigma_t)
$$
to $\R$.

If instead $\pi$ has a GT-reducible entry $x$ then we again remove $\pi$ from $P$, but now we search through our set of succession rules $\R$ and replace each instance of $(\pi)$ by the label we have given to the node $\st(\pi-x)$ (this label might not be $(\st(\pi-x))$ because $\st(\pi-x)$ may also have a GT-reducible  entry).  In other words, whenever a node labeled by $(\pi)$ would have been produced, now a node labeled by the same label as $\st(\pi-x)$ will be produced.  This does not change the isomorphism type of the tree because 
$$
T(B;\pi)=W(B;|\pi|,\pi)\cong W(B;|\pi|,\pi-x)\cong T(B;\st(\pi-x)).
$$

We repeat this process until $P=\emptyset$.  If we ever reach this state then we know that the generating tree we have produced is isomorphic to $T(B)$.  We will prove shortly (Theorem~\ref{finlabel}) that we do reach this state when $B$ contains both a child of an increasing permutation and a child of a decreasing permutation.

Before that, let us illustrate the process with the tree $T(132,3241)$.  We start with $P=\{1\}$ and $\R=\emptyset$.  Then we choose $1$ from $P$ and note that it does not have a GT-reducible entry (by our convention), so we remove $1$ from $P$ and add $12$ and $21$ to $P$, giving us $P=\{12,21\}$.  We also add the rule
$$
(1)\leadsto(12)(21)
$$
to our set of rules $\R$.

Now choose $12$ from $P$.  First we check to see if $1$ is a GT-reducible entry.  Since $\Binf=4$, Proposition~\ref{reinsert} shows that we only need to test whether $\ell_r(W(\{132,3241\};2,2))$ and $\ell_r(T(\{132,3241\};12))$ agree for $r=1,2,3$.  The number are shown in the following chart.
$$
\begin{array}{c|c|c}
r&\ell_r(W(\{132,3241\};2,2))&\ell_r(T(\{132,3241\};12))\\
\hline\hline
1&2&2\\
2&5&5\\
3&13&13
\end{array}
$$
From this chart we may conclude that $1$ is a GT-reducible entry, so we remove $12$ from $P$ and we replace $(12)$ by $(1)$ in all of our rules.  After this we have $P=\{21\}$ and $\R$ contains the single rule
$$
(1)\leadsto(1)(21).
$$

We must now choose $21$ from $P$.  However, $21$ does not have any removable entries, so it is not GT-reducible.  We therefore remove $21$ from $P$, add its children ($321$, $231$, and $213$) to $P$, and add the rule
$$
(21)\leadsto(321)(231)(213)
$$
to $\R$.

Using the same process as with $12$ we can find that both $321$ and $231$ have GT-reducible entries (the entry $2$ is GT-reducible for both permutations), so we replace $(321)$ and $(231)$ in all our rules by $(21)$.

At this point we have $P=\{213\}$.  This permutation also has a GT-reducible entry (again, $2$) so if we were to follow the pattern of the previous cases we would simply replace all instances of $(213)$ in our set of rules by $(12)$.  However, since $12$ is GT-reducible itself we instead replace instances of $(213)$ by $(1)$.

We are now done because $P=\emptyset$, and we have found the generating tree
$$
\begin{array}{llcl}
\mbox{Root:}	&	(1)&&\\
\mbox{Rule:}	&	(1)&\leadsto &(1)(21)\\
&(21)&\leadsto &(21)(21)(1).
\end{array}
$$
Up to relabeling, this is the same tree we found in Section~\ref{finlabel-intro}.  Via the transfer matrix method, one can use these rules to compute the generating function for $\Av(132,3241)$:
$$
\sum_{n\ge 1} |\Av_n(132,3241)|x^n
=
\frac{x(1-x)}{x^2-3x+1},
$$
and thus the sequence $|\Av_n(132,3241)|$ contains every second Fibonacci number.

\section{Proof of the main result}\label{finlabel-proof}

It remains to show that this procedure terminates.  In fact, if the nodes of $T(B)$ have arbitrarily large degrees, then we will never reach a state where $P=\emptyset$, because in that case either $12\cdots n$ or $n\cdots 21$ will never have a removable entry, let alone a GT-reducible entry.  Our central result, Theorem~\ref{finlabel} below, shows that if $B$ is finite and $T(B)$ has bounded degrees then this procedure will terminate.

We begin with a technical lemma that will be used to construct large sets of entries satisfying a condition stronger than removability.

\begin{lemma}\label{manyremove}
Suppose that $B$ is a finite set of patterns containing both a child of an increasing permutation and a child of a decreasing permutation and fix a positive integer $r$.  Every sufficiently long $B$-avoiding permutation $\pi$ contains a set $X=\{x_1,\dots,x_r\}$ of $r$ distinct pairwise nonadjacent entries so that each $x_j$ is removable in the word $\pi-(X\setminus\{x_j\})$.
\end{lemma}
\begin{proof}
We prove the lemma by induction on $r$.  The base case $r=1$ is immediate because Theorem~\ref{bounded} implies that all sufficiently long permutations contain removable entries, and we may take $x_1$ to be any removable entry in $\pi$.

Let $\pi$ be a $B$-avoiding permutation of length $n$.  The $r=2$ case provides a nice illustration of our argument, so we examine it before moving on to the general case.  Let $y_1,\dots,y_s$ denote the removable entries in $\pi-x_1$ that are not adjacent to $x_1$ in $\pi$, and assume to the contrary that $x_1$ is not removable in $\pi-y_i$ for any $i\in[s]$.  Because $x_1$ is removable in $\pi$, at least one of the sites adjacent to $x_1$ in $\pi$ is inactive.  Now form the $B$-containing permutation $\sigma$ by inserting $n+1$ into $\pi$ in an inactive site adjacent to $x_1$.  Because $x_1$ is not removable in $\pi-y_i$, $\sigma-y_i$ is $B$-avoiding.  But this means that every copy of a pattern from $B$ in $\sigma$ must contain the entry $y_i$ and the entry $n+1$, so $s\le\Binf-1$.  On the other hand, Theorem~\ref{bounded} shows that we may take $s$ to be as large as we like so long as $n$ is sufficiently large, a contradiction.

Now suppose that $r$ is any integer at least $2$, and that we have found a set of entries $\{x_1,\dots,x_{r-1}\}$ satisfying the desired conditions.  We wish to find an entry $x_{r}$, not adjacent to any of the entries $x_1,\dots,x_{r-1}$, so that $\{x_1,\dots,x_{r}\}$ satisfies the desired conditions.  As in the $r=2$ case before, we begin by letting $y_1,\dots,y_s$ denote the removable entries of $\pi-x_1-\cdots-x_{r-1}$ that are not adjacent (in $\pi$) to any $x_i$.  Assume to the contrary that none of the $y_i$'s will function adequately as $x_{r}$, so for each $i\in[s]$ there is at least one $j\in[r-1]$ such that $x_j$ is not removable in $\pi-x_1-\cdots-x_{j-1}-x_{j+1}-\cdots-x_{r-1}-y_i$.  Choose one of these values to be denoted $j(i)$.

When we were trying to find an entry to serve as $x_2$ we built a single permutation $\sigma$.  This time we need to consider the $r-1$ permutations $\sigma_1,\dots,\sigma_{r-1}$ where $\sigma_j$ is formed by inserting $n+1$ into $\pi$ in an inactive site adjacent to $x_j$ and then removing $x_1,\dots,x_{j-1},x_{j+1},\dots,x_{r-1}$.  Each of these permtuations contains a pattern from $B$ and by our assumptions, $\sigma_{j(i)}-y_i$ avoids $B$ for each $i\in[s]$.  As before, we now ask how many different values of $i$ can share the same value $j(i)$.  The answer is the same: every copy of a pattern from $B$ in  $\sigma_{j(i)}$ must contain both $y_i$ and $n+1$, so at most $\Binf-1$ values of $i$ may share the same $j(i)$.  Therefore we need only have $s>(r-1)(\Binf-1)$, which we get if $n$ is sufficiently large, to guarantee that at least one of the $y_i$'s can serve as $x_{r}$, completing the proof of the claim.
\end{proof}

\begin{theorem}\label{finlabel}
Let $B$ be a finite set of patterns.  The pattern-avoidance tree $T(B)$ is isomorphic to a finitely labeled generating tree if and only if $B$ contains both a child of an increasing permutation and a child of a decreasing permutation.  Furthermore, if $T(B)$ satisfies these conditions then the algorithm presented in Section~\ref{finlabel-alg} will find a finitely labeled generating tree isomorphic to $T(B)$.
\end{theorem}
\begin{proof}
These conditions on $B$ are necessary by Theorem~\ref{bounded}.  To prove the other direction, it suffices to show that every sufficiently long permutation is GT-reducible.  Take $\pi$ to be a $B$-avoiding permutation of length $n$.  We prove the theorem by showing that if $r$ is sufficiently large then at least one of the removable entries $x_1,\dots,x_r$ guaranteed by Lemma~\ref{manyremove} must be GT-reducible.

If $x_j$ is left-removable in $\pi-(X\setminus\{x_j\})$, set $\epsilon_j=-$.  Otherwise $x_j$ must be right-removable in $\pi-(X\setminus\{x_j\})$ and we set $\epsilon_j=+$.  Suppose to the contrary that no $x_j$ is GT-reducible, and thus for every $j\in[r]$ there is some $v_j\in W(B;n,\pi-x_j)$ so that $w_j=\iota^{\epsilon_j}_{\pi,x_j}(v_j)$ contains a pattern from $B$.  In fact, Propositions~\ref{embed} and \ref{reinsert} show that we may assume $|w_j|\le n+\Binf-1$.

Each copy of a pattern from $B$ in $w_j$ must use the entry $x_j$ since $\partial_{x_j}(w_j)=v_j$ is $B$-avoiding.  Hence at most $\Binf$ different $x_j$'s may share the same $w_j$.  Therefore it would suffice to show that the number of $w_j$'s is bounded by some constant depending only on the set $B$.  To accomplish this we will show that each $w_j$ lies in the set
$$
\iota^{\epsilon_1}_{\pi,x_1}\circ\iota^{\epsilon_2}_{\pi-x_1,x_2}\circ\cdots\circ\iota^{\epsilon_r}_{\pi-x_1-\cdots-x_{r-1},x_r}(W(B;n,\pi-X)),
$$
by showing that for all $j\in[r]$,
$$
\iota^{\epsilon_j}_{\pi,x_j}
=
\iota^{\epsilon_1}_{\pi,x_1}\circ\iota^{\epsilon_2}_{\pi-x_1,x_2}\circ\cdots\circ\iota^{\epsilon_r}_{\pi-x_1-\cdots-x_{r-1},x_r}\circ\partial_{x_r}\circ\cdots\circ\partial_{x_{j+1}}\circ\partial_{x_{j-1}}\circ\cdots\circ\partial_{x_1}
$$
as maps from $W(B;n,\pi-x_j)$ to $W(\emptyset;n,\pi)$.  Proposition~\ref{embed} implies that this would follow from
\begin{eqnarray*}
\lefteqn{
\iota^{\epsilon_1}_{\pi,x_1}\circ\iota^{\epsilon_2}_{\pi-x_1,x_2}\circ\cdots\circ\iota^{\epsilon_r}_{\pi-x_1-\cdots-x_{r-1},x_r}=
}
\\ & &
\iota^{\epsilon_j}_{\pi,x_j}\circ
\iota^{\epsilon_1}_{\pi-x_j,x_1}\circ
\cdots\circ
\iota^{\epsilon_{j-1}}_{\pi-x_1-\cdots-x_{j-2}-x_j,x_{j-1}}\circ
\iota^{\epsilon_{j+1}}_{\pi-x_1-\cdots-x_j,x_{j+1}}\circ
\cdots\circ
\iota^{\epsilon_r}_{\pi-x_1-\cdots-x_{r-1},x_r},
\end{eqnarray*}
and this follows from Remark~\ref{commute} because the $x_i$'s are pairwise nonadjacent.

Theorem~\ref{bounded} gives a bound on the number of children a node in $T(B)$ may have.  Let $\Delta$ denote this bound.  Since $W(B;n,\pi-X)\cong T(B;\st(\pi-X))$, the number of nodes in $W(B;n,\pi-X)$ of height between $0$ and $\Binf-1$ is bounded by $1+\Delta+\Delta^2+\cdots+\Delta^{\Binf-1}$.  Therefore, since we have shown that each $w_j$ lies in the image of this set under the map $\iota^{\epsilon_1}_{\pi,x_1}\circ\iota^{\epsilon_2}_{\pi-x_1,x_2}\circ\cdots\circ\iota^{\epsilon_r}_{\pi-x_1-\cdots-x_{r-1},x_r}$, we have bounded the number of possible $w_j$'s, completing the proof.
\end{proof}

Theorem~\ref{finlabel} only applies when the set of forbidden patterns $B$ is finite.  It can be shown that this hypothesis is necessary.  First note that if $\beta_1,\beta_2\in B$ where $\beta_2$ contains a $\beta_1$ pattern then $T(B)=T(B\setminus\{\beta_2\})$, or in other words, the $\beta_2$ restriction is superfluous.  Thus we may always assume that $B=\min(B)$, where here $\min(B)$ denotes the minimal elements of $B$ with respect to the pattern containment ordering.  It is not immediately obvious that there are sets of permutations $B$ such that $\min(B)$ is {\it not} finite.  Equivalently: are there infinite antichains of permutations?  Indeed, as has been rediscovered numerous times during the past thirty-five years, there are.  The reader is referred to Atkinson, Murphy, and Ru\v{s}kuc~\cite{amr:pwocsop} for constructions and references to earlier work.

In particular, $\Av(321,4123)$ contains an infinite antichain.  Following the argument given in Murphy's thesis~\cite{maximillian} one can show that there is a set of the form $\Av(B)$ where $321,4123\in B$ that does not have a rational generating function.  In more detail, there is an infinite antichain $U\subseteq\Av(321,4123)$ with at most one member of each length.  If $U_1$ and $U_2$ are two different subsets of $U$ then $\Av(321,4123,U_1)$ and $\Av(321,4123,U_2)$ have different generating functions.  Because $U$ is infinite, this gives uncountably many distinct generating functions, and thus they cannot all be rational.  In particular, these trees cannot all be finitely labeled.

\section{Conclusion}\label{finlabel-appl}

The algorithm presented here is implemented in the Maple package {\sc FinLabel} available at \texttt{http://math.rutgers.edu/\~{}vatter/}.  This algorithm is --- up to symmetry --- only a special case of the algorithm of Zeilberger~\cite{z:wilf}.  More precisely, if {\sc FinLabel} can find a generating tree isomorphic to $T(B)$ then Zeilberger's package {\sc Wilf} can enumerate $B^{-1}$-avoiding permutations.  Our algorithm has the advantage that, in the cases that it can handle, it returns both a generating tree and generating function, whereas Zeilberger's algorithm only returns a polynomial time algorithm for computing $|\Av_n(B)|$.  On the other hand, Zeilberger's algorithm is applicable in many other situations.  For example, it can enumerate the $1234$-avoiding permutations, which are known to have a non-algebraic generating function (Gessel~\cite{gessel}).  Zeilberger's algorithm is extended in Vatter~\cite{wilfplus}.

Since the writing of this paper, the enumerative implication of Theorem~\ref{finlabel} has been generalized by Albert, Linton, and Ru\v{s}kuc~\cite{insertion}.  They introduce a correspondence between permutations and words called the insertion encoding and prove that if $B$ satisfies the hypotheses of Theorem~\ref{finlabel} (in fact, weaker hypotheses suffice for their theorem) then the set of insertion encodings of permutations from $\Av(B)$ forms a regular language.  It then follows from the theory of formal languages that $\Av(B)$ has a rational generating function.

We conclude with several results that can now be proved completely automatically using the {\sc FinLabel} package.  First we have results from the classical paper of Simion and Schmidt~\cite{ss:rp}:
\begin{footnotesize}
$$
\begin{array}{l|l}
B				&	\mbox{generating function for $\Av(B)$}\\[2.5pt]
\hline\hline
\{123, 213	\}		&	\frac{x}{1-2x}\rule{0pt}{15pt}\\[2.5pt]
\{123, 231	\}		&	\frac{-x(1-x+x^2)}{(x-1)^3}\\[2.5pt]
\{123, 321	\}		&	x+2x^2+4x^3+4x^4\\[2.5pt]
\{132, 231	\}		&	\frac{x}{1-2x}\\[2.5pt]
\{312, 231	\}		&	\frac{x}{1-2x}\\[2.5pt]
\{123,132,213\}		&	\frac{x(1+x)}{1-x-x^2}\\[2.5pt]
\{123,132,231\}		&	\frac{x}{(x-1)^2}\\[2.5pt]
\{123,132,321\}		&	x+2x^2+3x^3+x^4\\[2.5pt]
\{123,231,312\}		&	\frac{x}{(x-1)^2}\\[2.5pt]
\{132,213,231\}		&	\frac{x}{(x-1)^2}
\end{array}
$$
\end{footnotesize}

West~\cite{west:trees} undertook a systematic study of permutations that avoid one pattern of length three and another of length four.  The generating functions that our algorithm can rederive are listed in the chart below.
\begin{footnotesize}
$$
\begin{array}{l|l}
B				&	\mbox{generating function for $\Av(B)$}\\[2.5pt]
\hline\hline
\{123, 3214\}			&
	\frac{x(1-x)}{1-3x+x^2}\\[2.5pt]
\{123, 3241\}			&	
	\frac{x(1-3x+4x^2-x^3)}{(2x-1)(1-x)^3}\\[2.5pt]
\{123, 3421\}			&
	\frac{x(1-3x+5x^2-2x^3)}{(1-x)^5}\\[2.5pt]
\{123, 4321\}			&	
	x+2x^2+5x^3+13x^4+25x^5+25x^6\\[2.5pt]
\{132, 3214\}			&	
	\frac{x(1-2x+2x^2)}{1-4x+5x^2-3x^3}\\[2.5pt]
\{132, 3241\}			&	
	\frac{x(1-x)}{1-3x+x^2}\\[2.5pt]
\{132, 3421\}			&	
	\frac{x(1-3x+3x^2)}{(1-x)(2x-1)^2}\\[2.5pt]
\{132, 4321\}			&
	\frac{x(1-3x+5x^2-2x^3+x^4)}{(1-x)^5}\\[2.5pt]
\{213,1234\}			&
	\frac{x(1-x)}{1-3x+x^2}\\[2.5pt]
\{213,1243\}			&
	\frac{x(1-x)}{1-3x+x^2}\\[2.5pt]
\{213,1423\}			&
	\frac{x(1-x)}{1-3x+x^2}\\[2.5pt]
\{213,4123\}			&
	\frac{x(1-x)}{1-3x+x^2}
\end{array}
$$
\end{footnotesize}

Finally we have permutations that avoid two patterns of length four.  The following generating functions, recently computed by Kremer and Shiu~\cite{ks:len4}, can also be found using {\sc FinLabel}.
\begin{footnotesize}
$$
\begin{array}{l|l}
B				&	\mbox{generating function for $\Av(B)$}\\[2.5pt]
\hline\hline
\{1234, 3214\}		&
	\frac{x(1-3x)}{(x-1)(4x-1)}\\[2.5pt]
\{1234, 3241\}		&
	\frac{x(1-11x+54x^2-151x^3+268x^4-313x^5+234x^6-108x^7+29x^8-4x^9)}
	{(1-3x+x^2)(2x-1)^2(x-1)^6}\\[2.5pt]
\{1234, 3421\}		&
	\frac{x(1-7x+24x^2-44x^3+62x^4-39x^5+32x^6-19x^7+4x^8)}{(1-x)^9}\\[2.5pt]
\{1234, 4321\}		&
	x+2x^2+6x^3+22x^4+86x^5+306x^6+882x^7+1764x^8+1764x^9\\[2.5pt]
\{1243, 3214\}		&
	\frac{x(1-4x+5x^2-3x^3)(1-x)}{1-7x+17x^2-22x^3+13x^4-4x^5}\\[2.5pt]
\{1243, 3241\}		&
	\frac{x(1-9x+31x^2-49x^3+37x^4-14x^5+2x^6)}{(1-x)(1-4x+2x^2)(1-3x+x^2)^2}\\[2.5pt]
\{1243, 3421\}		&
	\frac{x(1-9x+34x^2-64x^3+64x^4-28x^5+4x^6)}{(x-1)(2x-1)^5}\\[2.5pt]
\{1423, 3214\}		&
	\frac{x(1-6x+12x^2-7x^3+2x^4)}{1-8x+22x^2-25x^3+10x^4-2x^5}\\[2.5pt]
\{1423, 3241\}		&
	\frac{x(2x-1)^2(1-x)}{1-7x+16x^2-16x^3+4x^4}\\[2.5pt]
\{3214, 4123\}		&
	\frac{x(1-3x)}{(x-1)(4x-1)}\\
\end{array}
$$
\end{footnotesize}

\bigskip
\noindent{\it Acknowledgment.\/}  I would like to thank Doron Zeilberger for encouraging and entertaining me throughout this endeavor, Michael Albert for suggestions that greatly improved the {\sc FinLabel} package and comments about the case where $B$ is infinite (for which thanks are also owed to Maximillian Murphy), and finally Bruce Sagan and the anonymous referees, whose suggestions significantly improved the presentation.

\bigskip
\bibliographystyle{acm}
\bibliography{finlabel}

\end{document}